\documentstyle[12pt,amscd]{amsart}
\newtheorem{thm}{Theorem}[section]
\newtheorem{cor}{Corollary}[section]
\newtheorem{lem}{Lemma}[section]
\newtheorem{prop}{Proposition}[section]

\theoremstyle{definition}
\newtheorem{defn}{Definition}[section]
\newtheorem{exmp}{Example}[section]

\theoremstyle{remark}
\newtheorem{claim}{Claim}[section]
\newtheorem{rem}{Remark}[section]

\begin{document}
\title{Free pencils on Divisors}

\author{Roberto Paoletti}
\address{Department of Mathematics \\
UCLA \\
Los Angeles, CA 90024}
\curraddr{Dipartimento di Matematica, Universit\'a di Pavia,
27100 Pavia, Italy}
\bigskip
\maketitle

%\input p1b.tex
%\input p1b0.tex
%\input p1b01.tex
%\input p1b1.tex
%\input p1b11.tex
%\input p1b12.tex
%\input p1bib.tex
%\end{document}

\section{\bf {Introduction}}

In algebraic geometry, it is rather typical that the
embedding of a variety $Y$ in another variety
$X$ forces strong constraints on the existence of free linear
series on $Y$.
For example, a classical result in plane curve theory states
that
the gonality of a smooth plane curve of degree $d$ is $d-1$
(\cite{acgh}).
It is then natural to
look for general statements of this flavor.

One particular case, which is quite well understood, is
the one where $Y$ is a divisor in $X$.
This problem has been studied by several researchers.
In particular,
a wide range of situations is dealt with by the following
result of Sommese (\cite{so:amp}):
\begin{thm}{(Sommese)} Let $Y\subset X$ be
an irreducible smooth ample divisor
and let $\phi :Y@>>>B$ be a morphism
onto another projective manifold.
If $\dim (Y)\ge \dim (B)+2$ then $\phi$ extends to a morphism
$\psi :X@>>>B$.
\end{thm}
Serrano
(\cite{se:ext}) then studied the case where $B$ is a smooth curve
and $\dim (X)=2$ or $3$. Namely, he proved
the two following theorems:
\begin{thm}{(Serrano)}
Let $C$ be an irreducible smooth curve
contained in a smooth surface $S$. Suppose that there exists
a morphism $\phi :C@>>>\bold P^1$ of degree $d$.
If $C^2>(d+1)^2$, then there exists
a morphism $\psi :S@>>>\bold P^1$ extending
$\phi$.
\end{thm}

\noindent
and
\begin{thm}{(Serrano)}
Let $X$ be a smooth projective threefold, and let $S\subset X$
be a smooth very ample surface.
Let $\phi :S@>>>\bold P^1$ be a morphism with connected
fibers.
Let $g(F)$ be the arithmetic genus of a fiber and
set $d=F\cdot S$.
If $S^3>(d+1)^2$ and
$dim H^0(X,\cal O_X(S))\ge 3d+3+2g(F)$,
then $\phi$ extends to a morphism $\psi :X@>>>\bold P^1$.
\end{thm}

Actually, Serrano proves more, in the sense that he shows how
these statements imply analogous ones with $\bold P^1$ replaced
by a general smooth curve $B$, and he can also replace
the above numerical conditions by weaker ones if $S+K_X$ is
a numerically even divisor.
His argument is based on Miyaoka's vanishing theorem combined with
a refinement of Bompieri's method.
Furthermore,
Serrano applies the above results and methods to the
study of the ampleness of the adjoint divisor.

\bigskip

On the other hand, in a celebrated theorem
Reider \cite{re:vbls} has shown how adjunction
problems on surfaces can be exaustively studied using vector
bundle methods. His argument is based on an application
of Bogomolov's instability theorem.
Furthermore,
Reider himself has also given a proof along these lines
of a statement close to Serrano's theorem for surfaces (\cite{re:app}).
Also in light of Serrano's result for threefolds,
it is therefore reasonable to expect that
methods
of this type should
be applicable to obtain some more general statement about
the extension of linear series on a divisor.
Our result in this direction is the following:

\begin{thm}
($char(k)=0$)
Let $X$ be a smooth projective $n$-fold,
and let $Y\subset X$ be a reduced irreducible divisor. If
$n\ge 3$ assume that $Y$ is ample, and if $n=2$ assume
that $Y^2>0$ (so that in particular it is at least nef).
Let $\phi :Y@>>>\bold P^1$ be a morphism, and let $F$ denote
the numerical class of a fiber.

\noindent
(i) If
$$F\cdot Y^{n-2}<\sqrt {Y^n}-1,$$
then there exists a morphism $\psi :X@>>>\bold P^1$
extending $\phi$. Furthermore, the restriction
$$H^0(X,\psi ^{*}\cal O_{\bold P^1}(1))@>>>
H^0(Y,\phi ^{*}\cal O_{\bold P^1}(1))$$
is injective. In particular, $\psi$ is linearly normal if
$\phi$ is.

\noindent
(ii) If
$$F\cdot Y^{n-2}=\sqrt {Y^n}-1$$
and $Y^n\neq 4$, then either there exists an
extension $\psi :X@>>>\bold P^1$ of $\phi$, or else
we can find an effective divisor
$D$ on $X$ such that
$(D\cdot Y^{n-1})^2=(D^2\cdot Y^{n-2})Y^n$ and
$D\cdot Y^{n-1}=\sqrt {Y^n}$,
and an inclusion
$$\phi ^{*}\cal O_{\bold P^1}(1)\subset \cal O_Y(D).$$
\end{thm}
When applied to $n=2$ and $n=3$, this gives the
above statements of Serrano.
However, the hypothesis are weaker, because we are
not requiring $Y$ to be smooth
and we don't need the assumption about
the number of sections of $\cal O_X(Y)$.
Besides, we don't require $Y$ to be very ample (unlike
Serrano's statement for $n=3$).
We have furthermore a
description of what happens in the boundary situation; for example,
$(d+1)^2=C^2$ is the case of a minimal pencil on a smooth plane curve
(of degree $d+1$).
With respect to Reider's result on surfaces,
the assumption that $Y^2\ge 19$ and that $Y$ be smooth is
not necessary. The above furthermore shows that
conclusion (a) in Proposition 2.15
of \cite{re:app} always occurs for $f<\sqrt {Y^2}-1$
(just take the Stein factorization of $X@>>>\bold P^1$), and
therefore the other possibilities
can only occur in the boundary case (ii).
This in turn gives more information about this
case.
For example, comparison with Reider's theorem shows that
when $f=\sqrt {Y^2}-1$,
under the additional hypothesis that the curve $Y$ be smooth and $Y^2\ge 19$,
if $\cal O_X(D)$ is not base point free then
it has exactly one base point.

As to $n\ge 4$, this is clearly weaker than Sommese's result except
that we are not requiring $Y$ to be smooth.

The argument provides a direct geometric construction of the extension,
as follows. If we let $A=\phi ^{*}\cal O_{\bold P^1}(1)$,
$V=\phi ^{*}H^0(\bold P^1,\cal O_{\bold P^1}(1))$, we can define
a rank two vector bundle $\cal F$ by the exactness of the sequence
$$0@>>>\cal F@>>>V\otimes \cal O_X@>>>A@>>>0.$$
In light of Bogomolov's instability theorem on a
$n$-dimensional variety,
the given numerical assumption implies that $\cal F$ is Bogomolov
unstable with respect to $Y$, and so we have a saturated destabilzing
line bundle $\cal L\subset \cal F$. Then $\cal L=\cal O_X(-D)$ for
some effective divisor on $X$, and hence we are reduced to arguing that
the numerology forces $D$ to move in a base point free pencil.

Using the relative version of the Harder-Narashiman filtration
(\cite{fl}) this concrete description can be adapted to families of
morphisms, and one can also prove a more general statement about
morphisms to arbitrary smooth curves.
\bigskip

Finally, using recent results of Moriwaki concerning
a version of the Bogomolov-Gieseker inequality
in prime charachteristic, the above statements can be generalized to
varieties defined
over a field of charachteristic $p$.
\bigskip

This paper covers part of the content of my Phd thesis at
UCLA.
I want to thank Robert Lazarsfeld, my advisor, for introducing me to
Algebraic Geometry and taking continuous interest in my progress.

I am also endebted to a number of people for valuable comments
and discussions; among them D. Gieseker, M. Green
and especially A. Moriwaki.

\section{\bf {Instability of rank two bundles}}

In this section we collect
some statements about instability of rank two vector
bundles on a smooth projective manifold.
References in this direction are,
for example, \cite{bo:st}, \cite{gi}
and \cite{mi:cc}. For the statements in charachteristic $p$,
we shall be using results from \cite{mo:fpb}.

Let us first assume $char(k)=0$. We shall keep this
convention until otherwise stated.
The basic result is given by the Bogomolov-Gieseker
inequality for semistable bundles:

\begin{thm}
Let $S$ be a smooth projective surface, and let $\cal
E$ be a rank two vector bundle on $X$ with Chern
classes $c_1(\cal E)$ and $c_2(\cal E)$.
If $H$ is any polarization and $\cal E$ is $H$-semistable,
then $c_1(\cal E)^2-4c_2(\cal E)\le 0$.
\end{thm}

\begin{defn} Let $X$ be any projective $n$-dimensional
manifold and let $\cal E$ be a rank two vector bundle
on $X$. Let $c_i(\cal E)\in A^i(X)$ be the Chern classes
of $\cal E$, $i=1$ and $2$.
Define the {\it discriminant of $\cal E$} as
$$\Delta (\cal E)=c_1(\cal E)^2-4c_2(\cal E)\in A^2(X).$$
\label{defn:discr}
\end{defn}

\begin{lem} Let $X$ be a smooth projective $n$-fold
and fix a polarization $H$ on $X$.
Consider a rank two vector bundle $\cal E$ on $X$
which is $H$-unstable. Suppose that $\cal L_1,\cal L_2
\subset \cal E$ are line bundles and set $e=deg_H(\cal E)$,
$l_1=deg_H(\cal L_1)$ and $l_2=deg_H(\cal L_2)$.
Suppose that $2l_i>e$, $i=1,2$ and that
$\cal L_2$ is saturated in $\cal E$.
Then $\cal L_1\subset \cal L_2$.
\end{lem}

{\it Proof.} Set $l=min\{l_1,l_2\}$. By assumption, we have
$2l>e$. Let
$$\cal Q=:\cal E/\cal L_2.$$
Then $\cal Q$ is a torsion free sheaf on $X$.
If $\cal L_1\not\subset \cal L_2$, then the induced morphism
$\cal L_1@>>>\cal Q$ is not identically zero, and therefore
it is generically nonzero.
This implies that the obvious morphism of vector bundles
$$\cal L_1\oplus \cal L_2@>>>\cal E$$
is generically surjective. Hence the line bundle
$\wedge ^2\cal E\otimes \cal L_1^{-1}\otimes \cal L_2^{-1}$
is effective, and therefore
$$e\ge l_1+l_2\ge 2l,$$
a contradiction.
$\sharp$
\bigskip

\begin{rem} The above argument still works if $2l_1\ge e$.
\end{rem}

\begin{cor} Let $\cal E$ be an $H$-unstable rank two vector
bundle on $X$. If $\cal L\subset \cal E$ is a saturated
destabilizing line bundle, then it is the maximal
$H$-destabilizing line bundle of $\cal E$.
$\cal L$ contains any $H$-destabilizing line bundle of $\cal E$.
\end{cor}

\begin{defn} Let $S$ be a smooth projective surface, and let
$N^1(S)$ be the vector space of all numerical equivalence classes of
divisors on $S$. The positive cone $K^{+}(S)\subset N^1(S)$
is described by the equations $D^2>0$ and $D\cdot H>0$ for some
(and hence for all) polarizations $H$ on $S$.
\end{defn}

If we apply this to the situation of Bogomolov's theorem, we have
\begin{cor} Let $S$ be a smooth projective surface and let $\cal
E$ be a rank two vector bundle on $S$ with $\Delta (\cal E)>0$.
Then there exists a sequence
$$0@>>>A@>>>\cal E@>>>\cal B\otimes \cal J_Z@>>>0,$$
where $Z\subset S$ is local complete intersection codimension
two subscheme and $A$ and $B$ are line bundles on $S$ such that
$A-B\in K^{+}(S)$. Furthermore, $A$ is the maximal destabilizing
line bundle of $\cal E$ with respect to any polarization on $X$.
\label{cor:devissage2}
\end{cor}
{\it Proof.} Fix any polarization $H$ on $S$. Since
$\cal E$ is $H$-unstable, there is an exact sequence
$$0@>>>A@>>>\cal E@>>>B\otimes \cal J_Z@>>>0$$
where $A$ is the maximal destabilizing subsheaf of $\cal E$,
and in particular it is saturated.
Hence we have $(A-B)\cdot H>0$. On the other hand, since
$c_1(\cal E)=A+B$ and $c_2(\cal E)=A\cdot B+[Z]$, we also have
$$0<\Delta (\cal E)=(A+B)^2-4A\cdot B-4deg([Z])\le (A-B)^2.$$
This implies $A-B\in K^{+}(S)$, and therefore
$A$ strictly destabilizes $\cal E$ with respect to any polarization
on $S$. On the other hand, being saturated, it then has to
be the maximal destabilizing subsheaf of $\cal E$ with respect to
any polarization on $X$.
$\sharp$
\bigskip

We now want to generalize the above results to higher dimensional
varieties.
We start by recalling the following fundamental result of
Mumford-Mehta-Ramanathan (cfr \cite{mi:cc}):
\begin{thm} Let $X$ be a smooth projective $n$-fold, and let
$H$ be a polarization on $X$.
Suppose that $\cal E$ is a vector bundle on $X$. If $m\gg 0$,
and $Y\in |mH|$ is general, then the maximal destabilizing
subsheaf of $\cal E|_Y$ is the restriction to $Y$ of the
maximal destabilizing subsheaf of $\cal E$.
\label{thm:mumera}
\end{thm}
\begin{rem} By the {\it maximal destabilizing subsheaf} of $\cal E$
one means the first term $\cal E_1$ of the Harder-Narashiman filtration
of $\cal E$.
If $\cal E$ is semistable,
$\cal E_1=\cal E$.
\end{rem}

We generalize definition 1.3 as follows:
\begin{defn} Let $X$ be a smooth projective $n$-fold, and let
$H$ be a polarization on $X$. Denote by $N^1(X)$ the vector space
of all numerical equivalence classes of divisors on $X$.
Then the {\it $H$-positive cone} $K^{+}(X,H)\subset N^1(X)$
is described by the equations $D^2\cdot H^{n-2}>0$ and $D\cdot
H^{n-1}>0$. Note that this implies $D\cdot H^{n-2}\cdot L>0$ for
any other polarization $L$ on $X$.
\end{defn}

We then have:
\begin{thm} Let $X$ be a smooth projective $n$-fold and let $H$
be a fixed polarization on $X$. Consider a rank two vector bundle $\cal E$
on $X$ of discriminant $\Delta (\cal E)$. If
$\Delta (\cal E)\cdot H^{n-2}>0$, then
there exists an exact sequence
$$0@>>>\cal A@>>>\cal E@>>>\cal B\otimes \cal J_W@>>>0$$
where $W\subset X$ is a (possibly empty) codimension two
local complete intersection subscheme, and $\cal A$ and $\cal B$ are
line bundles on $X$ such that $\cal A-\cal B\in K^{+}(X,H)$.
\label{thm:main}
\end{thm}
{\it Proof.} For $n=2$, this is
the content of Corollary \ref{cor:devissage2}.
For $n\ge 3$,
let $V\in |mH|$ be general, with $m\gg 0$.
We may assume that $V$ is a smooth irreducible surface,
and that the maximal $H$-destabilizing subsheaf of $\cal E|_S$
is the restriction of the maximal $H$-destabilizing subsheaf of
$\cal E$ (Theorem \ref{thm:mumera}). By the hypothesis,
$$\Delta (\cal E|_S)=\Delta (\cal E)\cdot mH>0.$$
Therefore, by induction $\cal E|_V$ is Bogomolov-unstable
with respect to $H|_V$,
and so there exists an exact sequence
$$0@>>> A@>>>\cal E@>>> B\otimes \cal J_Z@>>>0,$$
satisfying the conclusions of theorem 1.1.
Furthermore, by the above there is $\cal A\subset \cal E$
such that $\cal A|_S=A$. Being normal of rank one, $\cal A$ is
a line bundle.
$\sharp$
\bigskip

\begin{rem} Note the inequality $(\cal A-\cal B)^2
\cdot H^{n-2}\ge \Delta (\cal E)\cdot H^{n-2}$.
\label{rem:sat}
\end{rem}

\begin{defn}
Let $X$ be a smooth $n$-dimensional projective variety, and
let $H$ be an line bundle on $X$.
Consider a rank two vector bundle $\cal E$ on $X$.
We shall say that $\cal E$ is {\it Bogomolov-unstable} with respect
to $H$ if there exists a line bundle $\cal L\subset
\cal E$ such that $2c_1(\cal L)-c_1(\cal E)\in K^{+}(X,H)$.
Hence Theorem \ref{thm:main} can be rephrased by saying that if
$\Delta (\cal E)\cdot H^{n-2}>0$, then $\cal E$ is Bogomolov-unstable
with respect to $H$.
\label{defn:bogunst}
\end{defn}

Let us now come to the case of positive charachteristic.
The basic result is given here by Moriwaki's generalization
of the Bogomolov-Gieseker inequality (\cite{mo:fpb}).
Before stating his theorem, we need the following:

\begin{defn} Let $X$ be a smooth projective $n$-fold
and let $H$ be an ample line bundle on $X$.
Let $\cal E$ be a rank two vector bundle on $X$. We say that
$\cal E$ is weakly $\mu$-semistable
w.r.t. $H$ if for any proper subsheaf
$\cal F\subset \cal E$ there exists an ample divisor $D$ on $X$
such that
$\mu (\cal F,H,D)\le \mu (\cal E,H,D)$, where
for a sheaf $\cal G$ we set $\mu (\cal G,H,D)=:\dfrac {c_1(\cal G)
\cdot H^{n-2}\cdot D}
{rank(\cal G)}$.
\label{defn:mu}
\end{defn}

\begin{rem} In any charachteristic, if $\cal E$ is
Bogomolov-unstable
w.r.t. $H$ (definition \ref{defn:bogunst}), then it is not
$\mu$-semistable. On the other hand, if $\cal E$
is not $\mu$-semistable w.r.t. $H$ and $\Delta (\cal E)\cdot H^{n-2}
>0$, then it is necessarily Bogomolov-unstable w.r.t. $H$.
\label{rem:mub}
\end{rem}

\begin{defn}
Let $X$ be a smooth projective $n$-fold, and let
$H$ be an ample line bundle on $X$,
and let $Nef(X)\subset N^1(X)$ denote the nef cone of
$X$.
Set
$$\sigma (H)=inf_{D\in Nef(X)}\Big \{\frac{(K_X\cdot D\cdot H^{n-2})^2}
{D^2\cdot H^{n-2}}\Big \}.$$
We agree to take the above ratio equal to
$\infty$ when the denumerator vanishes.
\label{defn:sigma}
\end{defn}

\begin{thm} (Moriwaki) Let $X$ be a smooth projective $n$-fold
over an algebraically closed field of
charachteristic $p>0$.
Assume that $X$ is not uniruled. Let $H$ be a polarization on $X$,
and let $\cal E$ be a rank two vector bundle on $X$.
Suppose that for all $0\le i<r$ the Frobenius pull-back $\cal E^{(i)}$
of $\cal E$ is weakly $\mu$-semistable with respect to $H$.
Then we have
$$\Delta (\cal E)\cdot H^{n-2}\le \frac {\sigma (H)}{(p^r-1)^2}.$$
\label{thm:mor}
\end{thm}

Furthermore, Moriwaki proves the following powerful
restriction lemma:
\begin{lem} ($char (k)\ge 0$)
Let $X$ be a smooth projective $n$-fold, and let $H$ be a very
ample line bundle on $X$. Suppose that $\cal E$ is a rank two vector bundle
on $X$, which is weakly $\mu$-semistable w.r.t. $H$.
Then for a general $Y\in |H|$ the restriction $\cal E|_Y$
is weakly $\mu$-semistable w.r.t. $H|_Y$.
\label{lem:A2}
\end{lem}

\begin{defn} Let $X$ be a smooth projective $n$-fold,
and let $H$ be an ample line bundle on $X$.
Define
$$\beta (H)=inf_{D\in Nef(X)}\Big \{\frac
{(D\cdot (H+K_X)\cdot H^{n-2})^2}{D^2\cdot H^{n-2}}\Big \}.$$
\label{defn:beta}
\end{defn}

\begin{cor} Let $X$ be a smooth projective $n$-fold,
with $n\ge 3$ on an
algebraically closed field,
and let $H$ be a very ample line
bundle on $X$.
Suppose that the general $Y\in |H|$
is not uniruled.
Let $\cal E$ be a rank two vector bundle
on $X$ such that
$$\Delta (\cal E)\cdot H^{n-2}>\dfrac {\beta (H)}{(p-1)^2}.$$
Then $\cal E$ is Bogomolov-unstable with respect to $H$,
i.e. there exists an exact sequence
$$0@>>>\cal A@>>>\cal E@>>>\cal B\otimes \cal J_Z@>>>0,$$
where $\cal A$ and $\cal B$ are line bundles on $X$,
$Z\subset X$ is a codimension two local complete intersection
and $\cal A-\cal B\in K^{+}(X,H)$.
\label{cor:devissagep}
\end{cor}

\begin{rem}
Although this is an immediate application
of Moriwaki's theorem \ref{thm:mor}, it is phrased in a way
that makes
it applicable to uniruled varieties.
Furthermore, observe that if $k$ is an uncountable
algebraically closed field and $X$ is a smooth
non-uniruled projective variety over $k$ with a
very ample line bundle $H$ on it,
the general element of
$|H|$ is not uniruled either.
\label{claim:uniruled}
In fact, since $k$ is uncountable, a variety $X$ over $k$ is uniruled if and
only if through a general point of $X$ there passes a rational curve
(\cite{mm}). But a general point in a general divisor of a very
ample linear series is a general point of $X$.
\end{rem}

{\it Proof.}
By Lemma \ref{lem:A2}, it is sufficient to
show that for general $Y\in |H|$ the restriction $\cal E|_Y$ is
Bogomolov-unstable with respect to $H|_Y$.
It is easy to deduce this fact from theorem \ref{thm:mor} and
the definition of $\beta (H)$.
$\sharp$
\bigskip

\begin{cor} Let $\cal E$ be a rank two vector bundle
on $\bold P^r_k$, where $k$ is an algebraically closed field
of charachteristic $p$. If $\Delta (\cal E)>0$,
then $\cal E$ is unstable.
\end{cor}
{\it Proof.}
For $r=2$, this is well-known. For
$r\ge 3$,
we apply corollary
\ref{cor:devissagep} taking the very ample line bundle
in the statement to be $\cal O_{\bold P^3}(4)$, so that
$\beta (H)=0$.
We can also proceed inductively from the case $r=2$
by applying Lemma \ref{lem:A2}.
$\sharp$
\bigskip

\section{\bf {Extension Of Pencils}}\label{section:ext}

Let $Y\subset X$ be an inclusion of projective varieties,
and let $|L|$ be a base point free pencil on $Y$. It is natural to
look for conditions under which $|L|$
extends to $X$, in the spirit of the results of Sommese,
Serrano and Reider (\cite{re:app}, \cite{so:amp}, \cite{se:ext}).
Our main result is the following:
\begin{thm} ($char(k)=0$)
Let $X$ be a smooth projective $n$-fold,
$n\ge 2$, and let
$Y\subset X$ be a reduced irreducible divisor.
If $n\ge 3$, assume
that $Y$ is ample, and if $n=2$ that $Y^2>0$
(so that in particular it is nef).
Suppose given a morphism $\phi :Y@>>>\bold P^1$
and let $F$ denote the numerical class of a fiber of $\phi$.

(i) Suppose that
$$F\cdot Y^{n-2}<\sqrt {Y^n}-1.$$
Then there exists a morphism $\psi :X@>>>\bold P^1$ extending
$\phi$,
and such that
$$H^0(X,\psi ^{*}\cal O_{\bold P^1}(1))\hookrightarrow
H^0(Y,\phi ^{*}\cal O_{\bold P^1}(1)).$$
In particular,
if $\phi$ is linearly complete, then so is
$\psi$.

(ii) Suppose $F\cdot Y^{n-2}=\sqrt {Y^n}-1$ and
$Y^n\neq 4$. Then either $\phi ^{*}\cal
O_{\bold P^1}(1)$ extends to a
base-point free pencil on $X$, or else there
exists an effective divisor $D$ on $X$ such that

\noindent
(a) the following equalities hold:
$$(D^2\cdot Y^{n-2})Y^n=(D\cdot Y^{n-1})^2$$
and
$$D\cdot Y^{n-1}=\sqrt {Y^n}.$$

\noindent
(b) there is an inclusion
$$\phi ^{*} \cal O_{\bold P^1}(1)\subset \cal O_Y(D).$$
\label{thm:ext}
\end{thm}

\begin{rem} If $Y$ is ample, the equalities in (a) of (ii)
can be phrased as follows. If $S\subset X$ is a smooth
complete intersection of $n-2$ divisor equivalent to multiples
of $Y$, then
$$D-\frac 1{\sqrt {Y^n}}Y\in Ker\{N(X)@>>>N(S)\}.$$
If $n=2$, this is just saying that $D\equiv _n\frac 1{\sqrt {Y^n}}Y$.
\end{rem}

{\it Proof.} Set
\begin{equation}
A=:\phi ^{*}\cal O_{\bold P^1}(1)
\label{eq:A}
\end{equation}
and let
\begin{equation}
V\subset H^0(Y,A)
\label{eq:V}
\end{equation}
be the pencil associated to $\phi$, i.e.
$V=\phi ^{*}H^0(\bold P^1,\cal O_{\bold P^1}(1))$.
Define a sheaf
$\cal F$ on $X$ by the exactness of
the sequence
\begin{equation}
0@>>>\cal F@>>>V\otimes \cal O_X@>>>A@>>>0.
\label{eq:F}
\end{equation}
Then $\cal F$ is a rank two vector bundle with Chern classes
$c_1(\cal F)=-Y$ and
$c_2(\cal F)=[A]$, where
$Y$ denotes the divisor class on $Y$ of an element
of the pencil $|A|$.
In particular, $[A]$ is represented by a fiber $F$ of $\phi$.
Therefore the discriminant of $\cal F$ (definition \ref{defn:discr}) is
\begin{equation}
\Delta (\cal F)=Y^2-4[A]
\label{eq:Delta}
\end{equation}
and so
\begin{equation}
\Delta (\cal F)\cdot Y^{n-2}=Y^n-4F\cdot Y^{n-2}.
\label{eq:Delta1}
\end{equation}
It is easy to check that
\begin{equation}
\sqrt {Y^n}-1\le \frac {Y^n}4
\label{eqn:easy}
\end{equation}
and by assumption we then have in particular that
$\Delta (\cal F)\cdot Y^{n-2}>0$, and therefore
$\cal F$ is Bogomolov-unstable with respect to $Y$
(definition \ref{defn:bogunst}).
Hence there exists a saturated invertible subsheaf
$$\cal L\subset \cal F$$
which is the maximal destabilizing subsheaf of $\cal F$
with respect to $(Y,\cdots,Y,L)$, for any
ample divisor $L$ on $X$ (theorem \ref{thm:main}).
Since $\cal L\subset \cal F\subset \cal O_X^2$, we can write
$$\cal L=\cal O_X(-D)$$
for some effective divisor $D$ on $X$.
The instability condition then reads
\begin{equation}
(Y-2D)\cdot Y^{n-1}\ge 0,
\label{eq:inst}
\end{equation}
with strict inequality holding if $Y$ is ample.
Furthermore, using the fact that $\cal L$ is saturated one can see
that
\begin{equation}
(Y-2D)^2\cdot Y^{n-2}\ge \Delta (\cal F)\cdot Y^{n-2}
\label{eq:inst1}
\end{equation}
(see remark \ref{rem:sat}) and if we set
$f=:F\cdot Y^{n-2}$ this can be rewritten as
\begin{equation}
f\ge D\cdot Y^{n-1}-D^2\cdot Y^{n-2}.
\label{eq:inst2}
\end{equation}
By assumption, we have $f<\sqrt {Y^n}-1$
and together with (\ref{eq:inst2}) this gives
$$D^2\cdot Y^{n-2}-1> D\cdot Y^{n-1}-\sqrt {Y^n}.$$
Applying the Hodge Index Theorem, we then get
\begin{equation}
\frac {(D\cdot Y^{n-1})^2}{Y^n} -1
>D\cdot Y^{n-1}-\sqrt {Y^n}.
\label{eq:hit}
\end{equation}
\begin{claim} $\cal L$ is saturated in $\cal O_X^2$.
\end{claim}
{\it Proof.} If not, there would exist an inclusion
$\cal O_X(Y-D)\subset \cal O_X^2$
(here we use the fact that $Y$ is reduced and irreducible)
and therefore
we should have
$$(D-Y)\cdot Y^{n-1}\ge 0.$$
Together with (\ref{eq:inst}), this would imply
$Y^n\le 0$, a contradiction.
$\sharp$
\bigskip

Hence we have an exact sequence of the form
\begin{equation}
0@>>>\cal O_X(-D)@>>>\cal O_X^2@>>>\cal O_X(D)
\otimes \cal J_Z@>>>0,
\label{eq:sat}
\end{equation}
where $Z\subset X$ is a codimension two local complete intersection.
Computing $c_2(\cal O_X^2)=0$ from the above sequence
we then get
$$D^2=[Z]$$
(equivalently, one might just observe that $Z$ is the
complete intersection of the two sections of $\cal O_X(D)$
coming from the above sequence).
Therefore under the assumptions of the theorem
either $Z=\emptyset$, or else $D^2\cdot Y^{n-2}>0$.

\begin{lem} $Z=\emptyset$
\end{lem}

{\it Proof.} Suppose, otherwise, that
$D^2\cdot Y^{n-2}>0$.
In this case the Hodge Index Theorem yields
$$(D\cdot Y^{n-1})^2\ge
(D^2\cdot Y^{n-2})Y^n\ge Y^n$$
and therefore
$$D\cdot Y^{n-1}\ge \sqrt {Y^n}.$$
Therefore the right hand side of (\ref{eq:hit}) is nonnegative.
We can rewrite (\ref{eq:hit}) as
$$\frac {(D\cdot Y^{n-1})^2}{Y^n}-1>
D\cdot Y^{n-1}-\sqrt {Y^n}=Y^n\{\frac {D\cdot Y^{n-1}}{Y^n}-
\frac 1{\sqrt {Y^n}}\}.$$
Let us now make use of the destabilizing condition
$Y^n\ge 2D\cdot Y^{n-1}$: we obtain
$$ \frac {(D\cdot Y^{n-1})^2}{Y^n}-1>
2\frac {(D\cdot Y^{n-1})^2}{Y^n}-2\frac {D\cdot Y^{n-1}}{\sqrt Y^n}$$
and this leads to $0>\Big \{
\dfrac {D\cdot Y^{n-1}}{\sqrt {Y^n}}-1\Big \}^2$, absurd.
$\sharp$
\bigskip

Since $Z=\emptyset$, $\cal O_X(-D)@>>>\cal O_X^2$ never drops
rank, and therefore neither does $\cal O_X(-D)@>>>\cal F$. Hence we have
a commutative diagram
\begin{equation}
\CD
  @.   @.  0   @.    0 @.  @.   \\
@.    @.         @VVV     @VVV     @.   \\
0@>>>\cal O_X(-D)@>>>\cal F@>>>\cal O_X(D-Y)@>>>0  \\
@.   @|               @VVV      @VVV        @.    \\
0@>>>\cal O_X(-D)@>>>V\otimes \cal O_X@>>>\cal O_X(D)@>>>0 \\
@.    @.  @VVV     @VVV   @.   \\
 @.  @.   A @=  \cal O_Y(D)  @.  @.  \\
@.   @.   @VVV   @VVV     @.  \\
  @.   @.  0  @.  0 @. @.
\endCD
\label{eq:bigcd}
\end{equation}
from which we see that
$$A\simeq \cal O_Y(D).$$
Furthermore, since
$$(D-Y)\cdot Y^{n-1}\le 2D\cdot Y^{n-1}-Y^n<0$$
by the destabilizing condition (\ref{eq:inst}), we have
$H^0(X,\cal O_X(D-Y))=0$ and therefore an
injection
$$H^0(X,\cal O_X(D))\hookrightarrow H^0(Y,\cal O_Y(D-Y)).$$
Since $\cal O_X(D)$ is a quotient of
$\cal O_X^2$, it is globally generated and $V$ gives
a base point free pencil of sections of
$\cal O_X(D)$.
$\sharp$
\bigskip

{\it Proof of (ii)}
Suppose now that $F\cdot Y^{n-2}=\sqrt {Y^n}-1$. It is easy
to see that if $Y^n\neq 4$ then
the inequality (\ref{eqn:easy}) is strict, and therefore
$\cal F$ is still Bogomolov
unstable with respect to $Y$.
Arguing exactly as in the proof of the previous lemma
we get:
\begin{lem} Either $Z=\emptyset$, or else
the following equalities hold:
$$(D\cdot Y^{n-1})^2=(D^2\cdot Y^{n-2})Y^n$$
and
$$D\cdot Y^{n-1}=\sqrt {Y^n}.$$
\end{lem}
To complete the argument, observe that a variant of the
commutative diagram (\ref{eq:bigcd}) gives the exact sequence
$$0@>>>\cal O_X(D-Y)\otimes \cal J_W@>>>
\cal O_X(D)\otimes \cal J_Z@>>>\cal O_Y(D)\otimes \cal I_{Z\cap
Y}@>>>0$$
where $\cal I$ denotes an ideal sheaf on $Y$.
Therefore we also get an isomorphism
$A\simeq \cal O_Y(D)\otimes \cal I_{Z\cap Y}.$
$\sharp$
\bigskip

\begin{cor} (Serrano) Let $S$ be a smooth projective surface
and let $C\subset S$ be an irreducible smooth curve with $C^2>0$.
Then either
$$gon(C)\ge \sqrt {C^2}-1$$
or else for every minimal pencil $A$ on $C$ there exists
a base point free pencil $\cal O_S(D)$ on $S$ such that
$$A\simeq \cal O_C(D).$$
\end{cor}

\begin{cor} Let $C\subset
\bold P^2$ be a smooth curve of degree $d$. Then
$$gon(C)= d-1.$$
Furthermore, any base point free pencil on $C$ is
given by projecting through a point of $C$.
\label{cor:plane}
\end{cor}

{\it Proof.} The bound in the theorem gives
$gon(C)\ge d-1$. On the other hand projecting from a point of
$C$ shows that equality must hold.
Let $A$ be any minimal pencil on $C$. We may assume
that $d>2$.
We are then in the boundary situation
$f=\sqrt {C^2}-1$ (case (ii) of Theorem \ref{thm:ext}, $n=2$).
Hence we must have an inclusion
$A\subset \cal O_C(H)$, which shows that $A$ has the form
\begin{equation}
A=\cal O_C(H-P)
\label{eq:proj}
\end{equation}
for some $P\in C$.
But (\ref{eq:proj}) is saying exactly that
$A$ is the pull back of the hyperplane bundle on $\bold P^1$
under the morphism given by projection from $P$.
Hence all the minimal pencils are
obtained in this way.
$\sharp$
\bigskip

\begin{exmp} Let us apply the Theorem to the
gonality of Castelnuovo extremal curves
in $\bold P^3$.
If $C$ has even degree $d=2a$,
then $C$ is the complete intersection of a quadric $S$ and
an hypersurface of degree $a$.
Suppose that $S$ is smooth.
Then either $gon(C)\ge \sqrt {C\cdot _SC}-1=\sqrt 2a-1$,
or else a minimal pencil is induced by a base point free
pencil on $S$. $C$ on $S\simeq \bold P^1\times \bold P^1$
is a curve of type $(a,a)$, and restriction to it
of the two rulings gives two pencils of degree $a=\frac d2$,
which is the well-known answer. The argument is the same for
even degree.
\end{exmp}
\begin{exmp} For an example with $n=3$, let $S\subset
\bold P^3$ be a smooth surface
of degree $s$ containing a line
$L$, and let
$\phi :S@>>>\bold P^1$ be induced by projection from
$L$. Then a straighforward computation shows
that $f>s\sqrt s-1$.
\end{exmp}

We now give an application to
singular plane curves.

\begin{cor} Let $C\subset \bold P^2$
be a reduced irreducible curve of degree $d$, and
suppose that the only singularities of $C$ are
ordinary singular points $P_1,\cdots,P_k$
of multiplicities $m_1,\cdots,m_k$, respectively.
Let $m=max\{m_i\}$ and denote by $\tilde C$ the normalization
of $C$.
Suppose that
$d^2>\sum _im_i^2$. Then
$$gon(\tilde C)\ge
min\Big \{\sqrt {d^2-\sum _im_i^2}-1,d-\sqrt {\sum _im_i^2}
\Big \}.$$
\end{cor}

{\it Proof.} Let
$$f:X@>>>\bold P^2$$
be the blow up of $\bold P^2$
at $P_1,\cdots,P_k$,
$$E_i=f^{-1}P_i$$
for $i=1,\cdots,k$ be the
exceptional divisors,and let
$\tilde C\subset X$ be the proper
transform of $C$.
Then $\tilde C$ is an irreducible smooth curve and
$$\tilde C\in |dH-\sum _{i=1}^km_iE_i|.$$
Therefore we have
$$\tilde C^2=d^2-\sum _{i=1}^km_i^2>0$$
by assumption, and the hypothesis of the theorem
are satisfied.
Hence either
$$gon(\tilde C)\ge \sqrt {\tilde C^2}-1,$$
or else there exists an effective divisor $D$ on $X$ moving in a
base point free pencil and inducing a minimal pencil
on $\tilde C$.
We may then assume that $D$ has the form
$$D=xH-\sum _ia_iE_i$$
with $x>0$ and all the $a_i\ge 0$.
The condition $D^2=0$ then gives
$$x=\sqrt {\sum _ia_i^2}.$$
Hence
$$D\cdot \tilde C=xd-\sum _ia_im_i\ge xd-\sqrt {\sum _ia_i^2}\sqrt {
\sum _im_i^2}= $$
$$=xd-x\sqrt {\sum _im_i^2}\ge d-\sqrt {\sum _im_i^2}.$$
The statement follows.
$\sharp$
\bigskip

\begin{exmp} Let us consider for example the case
of a reduced irreducible plane curve $C\subset \bold P^2$
whose only singularities are nodes $P_1,\cdots,P_{\delta}$.
Suppose that
$$4\delta <d^2.$$
Then by the Corollary
$$gon (\tilde C)\ge min\{\sqrt {d^2-4\delta}-1,d-2\sqrt {\delta}\}.$$
For example, if we also assume that
$$\delta <d-2$$
then
$$\sqrt{d^2-4\delta}-1>d-3,$$
and for any effective divisor $D=xH-\sum _ia_iE_i$
with $D^2=0$ it
is easy to see that
$D\cdot \tilde C\ge d-2$.
Since projecting from a node gives a pencil of degree $d-2$,
we then have
$$gon(\tilde C)=d-2.$$
\end{exmp}

%\vfill
%\eject

We now show how theorem 2.1 applies to
families of morphisms.

\begin{prop} Let $X$ be a smooth projective $n$-fold
and $Y\subset X$ be a reduced irreducible divisor in $X$.
Suppose that $Y$ is ample when $n>2$ and that $Y^2>0$ when
$n=2$.
Let $\Phi :Y\times B@>>>\bold P^1$ be a family of morphisms
with $B$ smooth
and set $\phi _b=\Phi |_{Y\times \{b\}}$. Denote by $F$ the
numerical class of a fiber of $\phi _b$ (it is independent of
$b\in B$), and suppose that
$$F\cdot Y^{n-2}<\sqrt {Y^n}-1.$$
Then there exists a nonempty open subset
$T\subset B$ and a morphism
$$\Psi :X\times T@>>>\bold P^1$$
that restricts to $\Phi$ on $Y\times T$.
\label{prop:rel}
\end{prop}

{\it Proof.} Let
$$A=\Phi ^{*}\cal O_{\bold P^1}(1)$$
and
$$V=:\Phi ^{*}H^0(\bold P^1, \cal O_{\bold P^1}(1)).$$
Then we can define a rank two vector bundle on the smooth
variety $X\times B$ in the usual guise, by the exactness
of the sequence
\begin{equation}
0@>>>\cal F@>>>V\otimes \cal O_{X\times B}@>>>A@>>>0.
\label{eq:Frel}
\end{equation}
For $b\in B$ let us set $X_b=X\times \{b\}$
and $A_b=A|_{X_b}$.
Then we have $\cal F|_{X_b}\simeq \cal F|_b$,
where $\cal F_b=:Ker\{V\otimes \cal O_{X_b}@>>>A_b\}$.
Then $\cal F$ can be seen as a family of vector bundles
on $X$, with Chern classes $c_1(\cal F)=-Y$ and $c_2(\cal F)=
[A_b]$.
As in the proof of the Theorem, these vector bundles are
Bogomolov unstable with respect to $Y$.
Let $\cal L_b\subset \cal F_b$ be the maximal destabilizing
line bundle of $\cal F_b$. By the construction in Theorem
\ref{thm:ext}, the morphisms $\psi _b$ are associated to base point
free pencils of sections of $\cal L_b^{-1}$ induced by $V$.
Therefore, the proposition will follow once we show that the
line bundles $\cal L_b$ can be glued to a line bundle $\cal L
\subset \cal F|_{X\times T}$
on some open subset $X\times T$.
In fact, we have:
\begin{claim} For some nonempty open subset
$T\subset B$ there exists a line bundle $\cal L\subset
F|_{X\times T}$ such that $\cal L$ restricts to $\cal L_b$
on $X_b$, for each $b\in T$.
\end{claim}
{\it Proof} This follows from the relative
version of the Harder-Narashiman
filtration introduced in
\cite{fhs} and \cite{fl}.
$\sharp$
\bigskip

This proves the statement of the proposition.
$\sharp$
\bigskip

In his paper (\cite{se:ext}) Serrano expressed his results
about extensions in terms of morphisms to arbitrary smooth curves.
It seems in order to give here a corresponding
generalization of theorem 2.1.
\begin{defn} Let $B$ be a smooth curve. We shall denote
by $s(B)$ the smallest degree of a
nondegenerate plane birational model
of $B$, i.e. the smallest $k$ for which $B$ has a birational
$g^2_k$.
Nondegenerate only means that we agree to take
$s(\bold P^1)=2$.
\end{defn}

\begin{cor} Let $X$ and $Y$ satisfy the hypothesis of the
theorem, and let $\phi :Y@>>>B$ be a morphism to a
smooth curve. Denote by $F$ the numerical class of a fiber of
$\phi$, and suppose that
$$(s(B)-1)F\cdot Y^{n-2}<\sqrt {Y^n}-1.$$
Then there exists a morphism
$\psi :X@>>>B$ extending $\phi$.
\end{cor}
{\it Proof.} We adapt the argument in \cite{se:ext}, Lemma 3.2.
Let $f:B@>>>G\subset \bold P^2$ be a plane birational model
of $B$, of degree $s=s(B)$, and let $B^{*}\subset B$ be the inverse
image of the smooth locus of $G$.
For $b\in B^{*}$, let $\pi _b:B@>>>\bold P^1$ be the projection
from $f(b)$; $\pi _b$ is a morphism of degree $s-1$.
Consider the composition $\phi _b=\pi _b\circ \phi:Y@>>>\bold P^1$.
A fiber of $\phi _b$ is numerically equivalent to a sum
of $s-1$ fibers of $\phi$, and the numerical hypothesis then
imply, by theorem \ref{thm:ext}, that there exist extensions
$\psi _b:X@>>>\bold P^1$.
By Proposition \ref{prop:rel}, we can find a nonempty open subset
$T\subset B$ and a morphism
$$\Psi :X\times T@>>>\bold P^1$$
extending the morphism
$$\Phi :Y\times T@>>>\bold P^1$$
given by $\Phi (y,b)=\phi _b (y)$.
{}From this one sees that, if
$$X@>\gamma _b>>\Delta _b@>g_b>>\bold P^1$$
is the Stein factorization of $\phi _b$, then
$\Delta _b\simeq \Delta$
for some fixed curve $\Delta$ and all the morphisms
$\gamma _b$ can be identified.
Consider the morphism
$$h=(\gamma |_Y,\phi):Y@>>>\Delta \times \bold P^1.$$
It is easy to see that $\pi _1:h(Y)@>>>\Delta$ is an isomorphism.
Hence we can define
$\psi =\pi _2\circ \pi _1 ^{-1}\circ \gamma$.
$\sharp$

Let us consider now the case of prime charachteristic.
We give the corresponding version of Theorem \ref{thm:ext}.

\begin{thm} Let $k$ be an algebraically closed field
of charachteristic $p$, and let
$X$ be a smooth projective n-fold over
$k$. Let $Y\subset X$ be a reduced irreducible divisor,
and suppose that there
exists a morphism $\phi :Y@>>>\bold P^1$; let
$F$ denote the numerical class of a fiber of $\phi$.
Then $\phi$ can be extended to a morphism $\psi :X@>>>\bold P^1$
in the following situations:

(i)
$char(k)\neq 2,3$,
$n=2$, $Y^2>0$, $deg(F)<\sqrt {Y^2}-1$, $X$ not of general type.

(ii) $n=2$, $Y^2>0$, $X$ is not uniruled and
$$deg(F)<min \Big \{\sqrt {Y^2}-1,\frac 14Y^2-\frac 1{(p-1)^2}
\sigma_S\Big \}.$$

(iii) $n\ge 3$, $X$ is not uniruled and there exists a ample
line bundle
$H$ on $X$, such that $Y\equiv lH$
and
$$F\cdot Y^{n-2}<min\Big \{\sqrt {Y^n}-1,
\frac 14Y^n-\frac {l^{n-2}}{(p-1)^2}\sigma (H)\Big \}.$$

(iv) $n\ge 3$ and there exists a very ample line bundle $H$ on
$X$ such that $Y\equiv lH$, the general
$Z\in |H|$ is not uniruled and
$$F\cdot Y^{n-2}<min \Big \{\sqrt {Y^n}-1,
\frac 14Y^n-\frac {l^{n-2}}{(p-1)^2}
\beta (H)\Big \}.$$
\label{thm:extp}
\end{thm}

\begin{rem} The definitions of $\sigma (H)$ and $\beta (H)$
are given in section 2 (definitions \ref{defn:sigma} and
\ref{defn:beta}).
If $n=2$, $\sigma $ does not depend on $H$,
and we denote it by $\sigma _S$.
\end{rem}

{\it Proof.}
As to (i), that $\phi$ does not extend means that
$\Delta (\cal F)>0$ (cfr eq. (\ref{eq:F})),
but $\cal F$ is not Bogomolov-unstable. That this forces
$X$ to be of general type is the content of Theorem
7 of \cite{sb}.
For the other
statements, the argument is exactly the same as in the charachteristic
zero case, the extra assumptions being needed to apply
the results about unstable rank two bundles from section 2
(e.g., Corollary \ref{cor:devissagep}).
$\sharp$
\bigskip

\begin{rem} For $n=2$, we have to assume that $S$ is not
uniruled to apply the chrarachteristic $p$ version of
Bogomolov's theorem. However, in the case of $\bold P^2$
Bogomolov's theorem still holds (\cite{schw}). We can
therefore still argue as in Corollary \ref{cor:plane}
to deduce the classical statement about the gonality of plane curves.
\end{rem}

Theorem \ref{thm:extp} can be strengthened as follows
(see Theorem \ref{thm:mor}).
\begin{thm} Let notation be as in Theorem
\ref{thm:extp}, and suppose that
$F\cdot Y^{n-2}<\sqrt{Y^n}-1$. Assume that $X$ is not uniruled.
Let $r$ be the smallest positive integer such that
$$F\cdot Y^{n-2}<
min\Big \{\sqrt{Y^n}-1,
\frac {Y^n}4-\frac {l^{n-2}}{(p^r-1)^2}\sigma (H)\Big \}.$$
Then if $X^{\prime}\supset Y^{\prime}$ denote
the $(r-1)$-th Frobenius pull-backs of the varieties $X$ and $Y$,
there exists $\psi :X^{\prime}@>>>\bold P^1$ extending
the induced morphism $\phi ^{\prime}:Y^{\prime}
@>>>\bold P^1$.
\end{thm}

\end{document}